\documentclass{ifacconf}

\usepackage{graphicx}      
\usepackage{natbib}        
\usepackage{amsmath}
\usepackage{amssymb}

\begin{document}
\begin{frontmatter}

\title{On the port-Hamiltonian representation of systems described by partial differential equations}


\author[First]{M. Sch\"{o}berl} 
\author[Second]{A. Siuka} 

\address[First]{Johannes Kepler University Linz, Institute of Automatic Control and Control Systems Technology, Austria \\ (e-mail: markus.schoeberl@jku.at)}
\address[Second]{Johannes Kepler University Linz, Austria \\ (e-mail: a.siuka@liwest.at)}

\begin{abstract}                
We consider infinite dimensional port-Hamiltonian systems. Based on a power balance relation we introduce the port-Hamiltonian system
representation where we pay attention to two different scenarios, namely the non-differential operator case and the differential 
operator case regarding the structural mapping, the dissipation mapping and the in/output mapping. 
In contrast to the well-known representation on the basis of the underlying Stokes-Dirac structure 
our approach is not necessarily based on using energy-variables 
which leads to a different port-Hamiltonian representation of the analyzed partial differential equations.  
   \bigskip

\textit{NOTICE: this is the author’s version of a work that was accepted for publication in ifac-papersonline.net. Changes resulting from the publishing process, such as peer review, editing, corrections, structural formatting, and other quality control mechanisms may not be reflected in this document. Changes may have been made to this work since it was submitted for publication. A definitive version was subsequently published in ifac-papersonline.net, DOI 10.3182/20120829-3-IT-4022.00001}
\end{abstract}


\end{frontmatter}

\section{Introduction}

Modeling of physical systems described by partial differential equations
(pdes) in a port-Hamiltonian framework has been treated extensively
over the last years. However, there are still many open issues and
there exists no unique representation. An important contribution in
this field was the introduction of Stokes-Dirac structures which allows
to analyze field theories in a port-Hamiltonian framework and to exploit
this system representation for the controller design, see e.g. \cite{SchaftMaschke,Machelli2004I,Machelli2004II,MachelliSIAM,Maschke2005}.
Roughly speaking, the key property of the Stokes-Dirac structure is
to represent the power balance relation of physical systems in an
formal way. This is achieved by combining so-called flows and efforts
in the domain and on the boundary of a system. Using this approach
the proper choice of energy variables is crucial and the Hamiltonian
is also a functional depending on energy variables. 

A different kind of a port-Hamiltonian representation, not focusing
on energy variables, but also based on a power conservation law, was
proposed in \cite{SchlacherHam2008,SchoeberlMCMDSPiezo2008}. The
proposed Hamiltonian framework can be seen as an extension of \cite{Olver}
by incorporating dissipation and boundary ports. This formulation
has been exploited for control issues in \cite{SchoeberlCDC} and
\cite{SiukaActa}. 

In this contribution we will present an enhancement of \cite{SchlacherHam2008,SchoeberlMCMDSPiezo2008}
in such a way that beside the input map as in \cite{SchoeberlMCMDSPiezo2008}
also the structural mapping $\mathcal{J}$ and the dissipation mapping
$\mathcal{R}$ may involve differential operators. 

In the approach based on Stokes-Dirac structures it is quite natural
that these mappings involve differential operators, however due to
the choice of energy variables the Hamiltonian density does not depend
on derivative variables, in contrast to the approach presented in
this paper. For instance in mechanics, the choice of energy variables
suggests to use the strain $\epsilon$ as the dependent coordinate,
whereas also the displacement $u$ can be used. This latter choice
leads to derivative coordinates when the potential energy has to be
stated, since the energy is a function of the strain and $\epsilon=\partial_{X}u$
is met in the case of a one-dimensional domain. This has severe consequences
in the application of the variational derivative applied to the Hamiltonian
density, leading to a different port-Hamiltonian representation, and
has also an impact on the computation of structural invariants, see
also \cite{SchoeberlCDC}.

This paper is organized as follows. In Section 2 some notation is
presented and the geometric objects that play a fundamental role in
the paper are introduced. The third section is dealing with the representation
of port-Hamiltonian systems described by pdes, where we focus on two
cases: the non-differential operator case and the differential operator
case. Two specific applications, the vibrating string and a simple
model of magnetohydrodynamics are analyzed in the fourth section to
demonstrate how the introduced differential operators have an impact
on the power balance in practice. The paper closes with some concluding
remarks.

\section{Notation and Preliminaries}

In this paper we will apply differential geometric methods and we
will use a notation that is similar to the one in \cite{Giachetta}.
To keep the formulas short and readable we will use tensor notation
and especially Einsteins convention on sums. We use the standard symbol
$\otimes$ for the tensor product, $\wedge$ denotes the exterior
product (wedge product), $\mathrm{d}$ is the exterior derivative,
$\rfloor$ the natural contraction between tensor fields. By $\partial_{\alpha}^{B}$
are meant the partial derivatives with respect to coordinates with
the indices $_{B}^{\alpha}$. Furthermore $C^{\infty}(\cdot)$ denotes
the set of the smooth functions on the corresponding manifold. Moreover
we will not indicate the range of the used indices when they are clear
from the context. Additionally pull back bundles are not indicated
to avoid exaggerated notation.

In the forthcoming we will consider bundle structures in order to
be able to separate dependent and independent coordinates. Let us
consider the bundle $\mathcal{X}\rightarrow\mathcal{D},\,(X^{A},x^{\alpha})\rightarrow(X^{A})$
where $x$ are the dependent and $X$ the independent coordinates.
The first jet manifold $\mathcal{J}^{1}(\mathcal{X})$ can be introduced
possessing the coordinates $(X^{A},x^{\alpha},x_{A}^{\alpha})$, where
the capital Latin indices $A,B$ are used for the base manifold $\mathcal{D}$
(independent coordinates) and $x_{A}^{\alpha}$ denote derivative
coordinates of first order (derivatives of the dependent coordinates
with respect to the independent ones). The jet structure also induces
the so-called total derivative
\[
d_{A}=\partial_{A}+x_{A}^{\alpha}\partial_{\alpha}+x_{AB}^{\alpha}\partial_{\alpha}^{B}
\]
acting on elements including first order derivatives and $x_{AB}^{\alpha}$
correspond to derivative coordinates of second order living in $\mathcal{J}^{2}(\mathcal{X})$,
the second jet manifold. We will mainly focus on the first order case
in the sequel, however the jet-structure and the total derivatives
can easily be extended to higher order cases. 

We will treat so-called densities in the sequel (a quantity that can
be integrated), where we pay special attendance to densities of the
form $\mathfrak{F}=\mathcal{F}\mathrm{d}X$ with $\mathcal{F}\in C^{\infty}(\mathcal{J}^{1}(\mathcal{X}))$
where $\mathrm{d}X$ denotes the volume element on the manifold $\mathcal{D}$,
i.e. $\mathrm{d}X=\mathrm{d}X^{1}\wedge\ldots\mathrm{\wedge d}X^{d}$
with $\mathrm{dim}(\mathcal{D})=d$. Since $\mathcal{F}\in C^{\infty}(\mathcal{J}^{1}(\mathcal{X}))$
is met in our case we restrict ourselves to first order theory and
additionally we denote by $F=\int_{\mathcal{D}}\mathfrak{F}$ the
integrated quantity, where of course the map $x=\Phi(X)$ leading
to $x_{A}=\partial_{A}\Phi(X)$ has to be plugged in to be able to
evaluate the integral properly. 

Based on the bundle structure $\mathcal{X}\rightarrow\mathcal{D}$
one can introduce several tangent structures, where for our considerations
the vertical tangent bundle $\mathcal{V}(\mathcal{X})$ and the cotangent
bundles representing differential forms will be used. Of special importance
will be the space $\Lambda_{1}^{d}(\mathcal{X})$ induced by the bundle
structure $\mathcal{X}\rightarrow\mathcal{D}$ 
\begin{eqnarray*}
\Lambda_{1}^{d}(\mathcal{X}) & = & \mathcal{T}^{*}(\mathcal{X})\wedge(\overset{d}{\Lambda}\mathcal{T}^{*}(\mathcal{D}))
\end{eqnarray*}
with a typical element $\omega=\omega_{\alpha}\mathrm{d}x^{\alpha}\wedge\mathrm{d}X$.
It is worth stressing that the functions $\omega_{\alpha}$ may depend
on derivative coordinates, however as stated above to simplify the
notation no pull backs will be indicated in the forthcoming, i.e.
it should be clear from the context which order of derivative is included
in the differential form%
\footnote{To be more precise: If $\omega_{\alpha}\in\mathcal{J}^{1}(\mathcal{X})$
then $\omega\in(\pi_{0}^{1})^{*}\Lambda_{1}^{d}(\mathcal{X})$ with
$\pi_{0}^{1}:\mathcal{J}^{1}(\mathcal{X})\rightarrow\mathcal{X}$.
For simplicity we write $\omega\in\Lambda_{1}^{d}(\mathcal{X})$ and
suppress the pull back.%
}. An important object is the horizontal exterior derivative $\mathrm{d}_{h}$,
which meets $\mathrm{d}_{h}(\phi)=\mathrm{d}X^{A}\wedge d_{A}(\phi)$
acting on a differential form $\phi$, where $d_{A}(\phi)$ denotes
the Lie-derivative of $\phi$ with respect to $d_{A}$ (see the appendix
for more details concerning the relationship of $\mathrm{d}$ and
$\mathrm{d}_{h}$ and Stokes theorem in that context).

Furthermore we will treat linear differential operators (of order
$k$) that are of the following form
\[
\mathfrak{D}:\Lambda_{1}^{d}(\mathcal{X})\rightarrow\mathcal{V}(\mathcal{X})
\]
that maps an element $\Lambda_{1}^{d}(\mathcal{X})$ of jet-order
$p$ to an element $\mathcal{V}(\mathcal{X})$ of jet-order $p+k$.
In coordinates we have
\[
\mathfrak{D}(\omega)=\mathfrak{D}^{\alpha\beta\mathfrak{K}}d_{\mathfrak{K}}(\omega_{\alpha})\partial_{\beta}\,,\,\,\,\,\,\,\,\,\, d_{\mathfrak{K}}=d_{A_{k}}\circ\ldots\circ d_{A_{1}}.
\]
with $\omega\in\Lambda_{1}^{d}(\mathcal{X})$. The adjoint operator
$\mathfrak{D}^{*}$ follows by integration by parts and fulfills the
condition
\begin{equation}
\mathfrak{D}(\omega)\rfloor\varpi=\mathfrak{D}^{*}(\varpi)\rfloor\omega+\mathrm{d}_{h}(\mathfrak{d})\label{eq:dad}
\end{equation}
with $\omega,\varpi\in\Lambda_{1}^{d}(\mathcal{X})$, where $\mathfrak{d}$
is a bilinear expression involving the total derivatives up to order
$k-1$, see \cite{Olver}.

\section{System representation}

In this section we will introduce port-Hamiltonian systems described
by pdes based on a power balance relation. This means that the system
is introduced in such a way that the power balance relation together
with the structure of the equations represent the physical process.
We will consider two cases: (i) the Hamiltonian depends on derivative
coordinates, but the operators $\mathcal{J},\mathcal{R}$ and $\mathcal{G}$
(namely the interconnection, the dissipation and the input/output
maps) are just linear maps and (ii) where we relax the assumption
regarding $\mathcal{J},\mathcal{R}$ and $\mathcal{G}$ and consider
adequate differential operators $\mathfrak{J},\mathfrak{R}$ and $\mathfrak{G}$
instead. Before we are able to introduce the corresponding system
representations we need to introduce some geometric concepts which
will be exploited later on with respect of the derivation of the power
balance relation. 

Let us consider a vector field which is used to measure the change
of the density $\mathcal{F}$. We use a (generalized) vertical vector
field $v:\mathcal{X}\rightarrow\mathcal{V}(\mathcal{X})$ locally
given as $v=v^{\alpha}\partial_{\alpha}$ where $v^{\alpha}$ may
depend on derivative coordinates, together with its first jet-prolongation
$j^{1}(v)$ which reads as 
\begin{eqnarray*}
j^{1}(v) & = & v^{\alpha}\partial_{\alpha}+d_{A}(v^{\alpha})\partial_{\alpha}^{A}.
\end{eqnarray*}
Then we compute the Lie-derivative of the densitiy $\mathfrak{F}$
as it has been defined before with respect to the vector field $j^{1}(v)$
and we obtain the important relation 
\begin{eqnarray}
j^{1}(v)(\mathcal{F}\mathrm{d}X) & = & \left(v^{\alpha}(\partial_{\alpha}\mathcal{F}-d_{A}\partial_{\alpha}^{A}\mathcal{F})+d_{A}(v^{\alpha}\partial_{\alpha}^{A}\mathcal{F})\right)\mathrm{d}X\nonumber \\
 & = & \left(v^{\alpha}\delta_{\alpha}\mathcal{F}+d_{A}(v^{\alpha}\partial_{\alpha}^{A}\mathcal{F})\right)\mathrm{d}X.\label{eq:intbyparts}
\end{eqnarray}
Here the Euler Lagrange operator $\delta$ 
\begin{equation}
\delta:\mathcal{J}^{2}(\mathcal{X})\rightarrow\Lambda_{1}^{d}(\mathcal{X})\label{eq:EulerLag}
\end{equation}
comes into play, whose coordinate expressions involves the variational
derivative $\delta_{\alpha}$ which acts on $\mathcal{F}$ as
\[
\delta_{\alpha}\mathcal{\mathcal{F}}=\partial_{\alpha}\mathcal{F}-d_{A}\partial_{\alpha}^{A}\mathcal{F}.
\]

\begin{rem}
In fact (\ref{eq:EulerLag}) written out in coordinates is a map 
\begin{eqnarray*}
\delta_{\alpha}\mathcal{\mathcal{F}} & \rightarrow & \delta_{\alpha}\mathcal{\mathcal{F}}\mathrm{d}x^{\alpha}\wedge\mathrm{d}X
\end{eqnarray*}
with $\mathfrak{F}=\mathcal{F}dX$.
\end{rem}
Applying the Theorem of Stokes (see again the appendix) to (\ref{eq:intbyparts})
we find that
\begin{eqnarray}
\int_{\mathcal{D}}j^{1}(v)(\mathcal{F}\mathrm{d}X) & = & \int_{\mathcal{D}}v^{\alpha}\left(\delta_{\alpha}\mathcal{F}\right)\mathrm{d}X+\int_{\partial\mathcal{D}}v^{\alpha}\left(\partial_{\alpha}^{A}\mathcal{F}\right)\mathrm{d}X_{A}\nonumber \\
 & = & \int_{\mathcal{D}}v\rfloor\delta\mathfrak{F}+\int_{\partial\mathcal{D}}v\rfloor\delta^{\partial}\mathfrak{F}\label{eq:LieDiffPde}
\end{eqnarray}
with $\mathrm{d}X_{A}=\partial_{A}\rfloor\mathrm{d}X$ and the boundary
operator 
\begin{eqnarray*}
\delta^{\partial}\mathfrak{F} & = & \partial_{\alpha}^{A}\mathcal{F}\,\mathrm{d}x^{\alpha}\wedge\mathrm{d}X_{A}.
\end{eqnarray*}
The relation (\ref{eq:LieDiffPde}) will be of key interest in the
forthcoming, since it provides a natural decomposition of the expression
$\int_{\mathcal{D}}j^{1}(v)(\mathcal{F}\mathrm{d}X)$ into a term
on the domain $\mathcal{D}$ and one the boundary $\partial\mathcal{D}$. 

Important is furthermore the case when the generalized vector-field
$v$ is linked to the solution of a pde system (via its semi-group,
that $v$ may generate), then the formal change of $F=\int_{\mathcal{D}}\mathcal{F}\mathrm{d}X$
along solutions of a pde system can be computed as $\int_{\mathcal{D}}j^{1}(v)(\mathcal{F}\mathrm{d}X)$
(provided all operations are admissible), which will be denoted by
\begin{eqnarray*}
\dot{F} & = & \int_{\mathcal{D}}j^{1}(v)(\mathcal{F}\mathrm{d}X)
\end{eqnarray*}
in this special case.

\subsection{The non-differential Operator case:}

We introduce a port controlled Hamiltonian system on a bundle $\mathcal{X}\rightarrow\mathcal{D}$
in the non-differential operator case in the form of 
\begin{equation}
\begin{array}{ccl}
\dot{x} & = & (\mathcal{J}-\mathcal{R})(\mathrm{\delta}\mathfrak{H)}+u\rfloor\mathcal{G}\\
y & = & \mathcal{G}^{*}\rfloor\mathrm{\delta}\mathfrak{H}.
\end{array}\label{eq:HamTimeInvStandardCF-1}
\end{equation}
together with appropriate boundary conditions and with a first order
Hamiltonian $\mathfrak{H}=\mathcal{H}\mathrm{d}X$ with $\mathcal{H}\in C^{\infty}(\mathcal{J}^{1}(\mathcal{X}))$.
(Here it is worth mentioning that we restrict ourselves to the first
order case). 
\begin{rem}
It is worth stressing again, that of course $\delta$ is a differential
operator, and the phrase 'non-differential operator case' only refers
to the maps $\mathcal{J},\mathcal{R}$ and $\mathcal{G}$.
\end{rem}
We make use of the linear maps (over vector bundles) 
\begin{eqnarray*}
\mathcal{J},\mathcal{R}:\Lambda_{1}^{d} & \rightarrow & \mathcal{V}(\mathcal{X})
\end{eqnarray*}
where $\mathcal{J}$ is skew-symmetric and meets
\begin{eqnarray}
\mathcal{J}(\omega)\rfloor(\varpi)+\mathcal{J}(\varpi)\rfloor(\omega) & = & 0\label{eq:Jop}
\end{eqnarray}
and since $\mathcal{R}$ is symmetric and positive semi-definite we
have 
\begin{equation}
\mathcal{R}(\omega)\rfloor(\varpi)-\mathcal{R}(\varpi)\rfloor(\omega)=0\,,\,\,\,\,\,\,\mathcal{R}(\omega)\rfloor(\omega)\geq0\label{eq:Rop}
\end{equation}
with $\omega,\varpi\in\Lambda_{1}^{d}$. 

The input and the output bundles are given as $\mathcal{U}\rightarrow\mathcal{X}$
and $\mathcal{Y}\rightarrow\mathcal{X}$ and the map $\mathcal{G}$
and its adjoint $\mathcal{G}^{*}$ are of the form
\[
\mathcal{G}:\mathcal{U}\rightarrow\mathcal{V}\left(\mathcal{X}\right)\,,\qquad\mathcal{G}^{*}:\Lambda_{1}^{d}\rightarrow\mathcal{Y}.
\]
with 
\begin{equation}
(u\rfloor\mathcal{G})\rfloor\omega=u\rfloor(\mathcal{G}^{*}\rfloor\omega).\label{eq:Gop}
\end{equation}
Replacing $\mathfrak{F}$ by $\mathfrak{H}$ in (\ref{eq:LieDiffPde})
and setting $v=\dot{x}$ we obtain 
\begin{eqnarray}
\dot{H} & = & \int_{\mathcal{D}}(\mathcal{J}-\mathcal{R})(\mathrm{\delta}\mathfrak{H)}\rfloor\delta\mathfrak{H}+\int_{\mathcal{D}}(u\rfloor\mathcal{G})\rfloor\delta\mathfrak{H}+\int_{\partial\mathcal{D}}\dot{x}\rfloor\delta^{\partial}\mathfrak{H}\nonumber \\
 & = & -\int_{\mathcal{D}}\mathcal{R}(\mathrm{\delta}\mathfrak{H)}\rfloor\delta\mathfrak{H}+\int_{\mathcal{D}}u\rfloor y+\int_{\partial\mathcal{D}}\dot{x}\rfloor\delta^{\partial}\mathfrak{H}\label{eq:HdotPdeoO}
\end{eqnarray}
which reflects the power balance, since the total change of the functional
$H$ along solutions of (\ref{eq:HamTimeInvStandardCF-1}), is made
up of dissipation, collocation on the domain and an expression corresponding
to a boundary port (if it exists) depending on the boundary conditions,
see \cite{SchoeberlMCMDSPiezo2008,Siukaphd}.
\begin{rem}
Since $\int_{\mathcal{D}}u\rfloor y$ in coordinates reads as $\int_{\mathcal{D}}u^{i}y_{i}\mathrm{d}X$
we interpret the outputs as density valued such that the fibers of
the vector bundle $\mathcal{U\rightarrow X}$ possess the base $e_{i}$
and the fibers of the vector bundle $\mathcal{Y\rightarrow X}$ possess
by duality the base $e^{i}\otimes\mathrm{d}X$. Similarly, this is
reflected also in the case of the pairing of $\mathcal{V}(\mathcal{X})$
with $\Lambda_{1}^{d}$ which also takes its values in the densities
since with $v\in\mathcal{V}(\mathcal{X})$ and $\omega\in\Lambda_{1}^{d}$
we have $v\rfloor\omega=v^{\alpha}\omega_{\alpha}\mathrm{d}X$.
\end{rem}

\subsubsection{Casimir functionals}

Let us consider how to obtain so-called Casimir densities (or functionals)
where we restrict ourselves to the first order case only, i.e. $\mathfrak{C}=\mathcal{C}\mathrm{d}X$
with $\mathcal{C}\in C^{\infty}(\mathcal{J}^{1}(\mathcal{X}))$. Using
(\ref{eq:LieDiffPde}) where we replace $\mathcal{F}$ by $\mathcal{C}$,
setting $v=\dot{x}$ and $C=\int_{\mathcal{D}}\mathfrak{C}$ the relation
\[
\dot{C}=\int_{\mathcal{D}}(\mathcal{J}-\mathcal{R})(\mathrm{\delta}\mathfrak{H})\rfloor\delta\mathfrak{C}+\int_{\mathcal{D}}(u\rfloor\mathcal{G})\rfloor\delta\mathfrak{C}+\int_{\partial\mathcal{D}}\dot{x}\rfloor\delta^{\partial}\mathfrak{C}
\]
is obtained and a Casimir functional meets
\[
\dot{C}=\int_{\mathcal{D}}(u\rfloor\mathcal{G})\rfloor\delta\mathfrak{C}.
\]
This leads to the following two conditions for the Casimir density
\begin{eqnarray}
\delta_{\alpha}\mathcal{C}(\mathcal{J}^{\alpha\beta}-\mathcal{R}^{\alpha\beta}) & = & 0\label{eq:CasPde1}\\
\left.(\dot{x}^{\alpha}\partial_{\alpha}^{A}\mathcal{C})\right|_{\partial\mathcal{D}} & = & 0\label{eq:CasPde2}
\end{eqnarray}
which have to be fulfilled (independently of the Hamiltonian density
$\mathfrak{H}$). If in addition 
\begin{eqnarray*}
\int_{\mathcal{D}}(u\rfloor\mathcal{G})\rfloor\delta\mathbb{\mathfrak{C}} & = & 0
\end{eqnarray*}
is met, then the density is a conserved quantity. 

The concept of control by interconnection based on Casimir functionals
using the presented system representation has been applied to a heavy
chain system in \cite{SchoeberlCDC} and to the Timoshenko beam in
\cite{SiukaActa}. The generation of Casimir functionals based on
the Stokes-Dirac structure and its consequences for control have been
treated for instance in \cite{Machelli2004II,Pasu}.

\subsection{The differential operator case:}

A port-Hamiltonian system with respect to the differential operator
case reads as
\begin{eqnarray}
\dot{x} & = & \left(\mathfrak{J}-\mathfrak{R}\right)\left(\delta\mathfrak{H}\right)+\mathfrak{G}\left(u\right)\label{eq:iPCHD_operator_case}\\
y & = & \mathfrak{G}^{*}\left(\delta\mathfrak{H}\right)\nonumber 
\end{eqnarray}
together with appropriate boundary conditions, see also \cite{Siukaphd}.

For this case the operators $\mathfrak{J}$,$\mathfrak{R:}\Lambda_{1}^{d}\rightarrow\mathcal{V}(\mathcal{X})$
are $r$-order linear vector valued differential operators where $\mathfrak{J}$
is a skew-adjoint operator according to
\begin{equation}
\mathfrak{J}\left(\omega\right)\rfloor\varpi+\mathfrak{J}\left(\varpi\right)\rfloor\omega=\mathrm{d}_{h}\left(\mathfrak{\mathfrak{j}}\right)\label{eq:J_operator}
\end{equation}
with 
\begin{eqnarray*}
\mathfrak{\mathfrak{j}} & = & \mathfrak{\mathfrak{j}}^{A}\mathrm{d}X_{A}\,,\,\,\,\,\,\omega,\varpi\in\Lambda_{1}^{d},
\end{eqnarray*}
cf. (\ref{eq:dad}). Furthermore, $\mathfrak{R}$ is a non-negative
self-adjoint operator, i.e.,
\begin{equation}
\mathfrak{R}\left(\omega\right)\rfloor\varpi-\mathfrak{R}\left(\varpi\right)\rfloor\omega=\mathrm{d}_{h}\left(\mathfrak{r}\right)\,,\qquad\mathfrak{R}\left(\omega\right)\rfloor\omega\geq0\,.\label{eq:R_operator}
\end{equation}
with $\mathfrak{r}=\mathfrak{r}^{A}\mathrm{d}X_{A}$. 

The input operator $\mathfrak{G}$ and its adjoint operator $\mathfrak{G}^{*}$
are maps according to
\begin{equation}
\mathfrak{G}:\mathcal{U}\rightarrow\mathcal{V}\left(\mathcal{X}\right)\,,\qquad\mathfrak{G}^{*}:\Lambda_{1}^{d}\rightarrow\mathcal{Y}\label{eq:G_diff_case}
\end{equation}
and they are linear $r$-order differential operators with respect
to the relation
\begin{equation}
\mathfrak{G}\left(u\right)\rfloor\omega=u\rfloor\mathfrak{G}^{*}\left(\omega\right)+\mathrm{d}_{h}\left(\mathfrak{g}\right)\,,\qquad\mathfrak{g}=\mathfrak{g}^{A}\mathrm{d}X_{A}\,,\label{eq:G_operator}
\end{equation}
Consequently, from (\ref{eq:G_operator}) we are able to derive the
relation

\begin{eqnarray}
\mathfrak{G}\left(u\right)\rfloor\delta\mathfrak{H} & = & u\rfloor\mathfrak{G}^{*}\left(\delta\mathfrak{H}\right)+\mathrm{d}_{h}\left(\mathfrak{g}\right)\nonumber \\
 & = & u\rfloor y+\mathrm{d}_{h}\left(\mathfrak{g}\right)\label{eq:GuyBound}
\end{eqnarray}
where the pairing $u\rfloor y$ will correspond to the domain port.
\begin{rem}
It should be noted that the relations (\ref{eq:J_operator}), (\ref{eq:R_operator})
and (\ref{eq:G_operator}) are the generalizations of (\ref{eq:Jop}),
(\ref{eq:Rop}) and (\ref{eq:Gop}) to the differential operator case,
where especially the horizontal derivative will lead to appropriate
boundary terms when the adjoint operator is taken into account.
\end{rem}
The formal change of the Hamiltonian functional $H$ along (\ref{eq:iPCHD_operator_case})
now takes the form
\begin{equation}
\dot{H}=\int_{\mathcal{D}}(\mathfrak{J}-\mathfrak{R})(\mathrm{\delta}\mathfrak{H)}\rfloor\delta\mathfrak{H}+\int_{\mathcal{D}}\mathfrak{G}\left(u\right)\rfloor\delta\mathfrak{H}+\int_{\partial\mathcal{D}}\dot{x}\rfloor\delta^{\partial}\mathfrak{H}\label{eq:HdotGeneralCase}
\end{equation}
and due to the involved differential operators it cannot be concluded
in a unique fashion how the power is extracted over the domain and/or
the boundary. It is easily checked for instance that plugging in (\ref{eq:GuyBound})
into (\ref{eq:HdotGeneralCase}) leads to collocation on the domain
and an additional boundary expression. These non-uniqueness depends
on the properties of the used differential operators in $\mathfrak{J},\mathfrak{R}$
and $\mathfrak{G}$ which also has an impact if Casimir functionals
are to be considered in this differential-operator case. To be able
to draw some conclusions also in that case we will restrict ourselves
to specific operators described in the following section:

\subsubsection{Specific Operators}

Motivated by the forthcoming applications we will consider two types
of operators. 

We introduce a second-order non-negative self-adjoint operator $\mathfrak{R}$
locally given as
\begin{equation}
\mathfrak{R}\left(\omega\right)=d_{A}(\mathfrak{R}^{\alpha\beta AB}\, d_{B}(\omega_{\beta}))\partial_{\alpha},\label{eq:R_2ndorder_local}
\end{equation}
with 
\begin{eqnarray*}
\mathfrak{R}^{\alpha\beta AB} & = & \mathfrak{R}^{\beta\alpha BA}\,,\,\,\,\,\,\,\,\mathfrak{R}^{\alpha\beta AB}\in C^{\infty}\left(\mathcal{X}\right).
\end{eqnarray*}
It satisfies the relation
\begin{equation}
\mathfrak{r}^{A}\mathrm{d}X_{A}=\bar{\mathfrak{R}}\left(\omega\right)\rfloor\varpi_{\partial}-\bar{\mathfrak{R}}\left(\varpi\right)\rfloor\omega_{\partial},\label{eq:R_2ndorder_operator}
\end{equation}
cf. (\ref{eq:R_operator}) with 
\begin{eqnarray*}
\omega_{\partial} & = & -\partial_{A}\rfloor\omega\,,\,\,\,\,\,\,\varpi_{\partial}=-\partial_{A}\rfloor\varpi
\end{eqnarray*}
where we used 
\begin{eqnarray*}
\bar{\mathfrak{R}}\left(\omega\right) & = & \mathfrak{R}^{\alpha\beta AB}\, d_{B}\left(\omega_{\beta}\right)\partial_{\alpha}.
\end{eqnarray*}

\begin{rem}
The condition $\mathfrak{R}^{\alpha\beta AB}\in C^{\infty}\left(\mathcal{X}\right)$
might be relaxed, and $C^{\infty}\left(\mathcal{X}\right)$ can be
replaced by an appropriate jet-space. Since, in our concrete application
this is not the case, we will not focus on that fact more detailed.
\end{rem}
The operator $\mathfrak{R}$ is non-negative since $\mathfrak{R}\left(\omega\right)\rfloor\omega$
can be written as 
\begin{equation}
-d_{A}\left(\omega_{\alpha}\right)\mathfrak{R}^{\alpha\beta AB}d_{B}\left(\omega_{\beta}\right)\mathrm{d}X+\mathrm{d}_{h}\left(\mathfrak{B_{R}}\right)\label{eq:Rintbyparts}
\end{equation}
with the boundary expression 
\begin{eqnarray*}
\mathfrak{B_{R}} & = & \omega_{\alpha}\mathfrak{R}^{\alpha\beta AB}\, d_{B}(\omega_{\beta})\mathrm{d}X_{A}.
\end{eqnarray*}
The non-negativity of the operator follows if
\begin{equation}
-d_{A}\left(\omega_{\alpha}\right)\mathfrak{R}^{\alpha\beta AB}d_{B}\left(\omega_{\beta}\right)\geq0\label{eq:non-negativity}
\end{equation}
is met.

Furthermore, we define a first-order input operator $\mathfrak{G}$
with
\begin{equation}
\mathfrak{G}\left(u\right)=\mathfrak{G}_{i}^{\alpha A}d_{A}\left(u^{i}\right)\partial_{\alpha}\label{eq:G_1storder_local}
\end{equation}
and, $\mathfrak{G}_{i}^{\alpha A}\in C^{\infty}(\mathcal{J}^{1}(\mathcal{X}))$
which meets
\begin{equation}
\mathfrak{G}\left(u\right)\rfloor\omega=u\rfloor\mathfrak{G}^{*}\left(\omega\right)+\mathrm{d}_{h}\left(\bar{\mathfrak{G}}\left(u\right)\rfloor\omega_{\partial}\right)\label{eq:G_1storder_operator}
\end{equation}
with
\[
\bar{\mathfrak{G}}\left(u\right)=\mathfrak{G}_{i}^{\alpha A}u^{i}\partial_{\alpha}\,,\quad\omega_{\partial}=-\partial_{A}\rfloor\omega.
\]
The relation (\ref{eq:G_1storder_operator}) takes the form of
\begin{eqnarray}
\omega_{\alpha}\mathfrak{G}_{i}^{\alpha A}d_{A}(u^{i})\mathrm{d}X & = & -u^{i}d_{A}(\omega_{\alpha}\mathfrak{G}_{i}^{\alpha A})\mathrm{d}X+\mathrm{d}_{h}(\mathfrak{B_{G}})\label{eq:Gintbyparts}
\end{eqnarray}
with 
\begin{eqnarray*}
\mathfrak{B_{G}} & = & \omega_{\alpha}\mathfrak{G}_{i}^{\alpha A}u^{i}\mathrm{d}X_{A}.
\end{eqnarray*}

\begin{rem}
As above the condition $\mathfrak{G}_{i}^{\alpha A}\in\mathcal{J}^{1}(\mathcal{X})$
can be relaxed in the same spirit as for the dissipation operator.
\end{rem}

\section{Examples}

In this section we will discuss the presented ideas based on physical
examples. We will present the vibrating string example in two fashions
-with and without structural damping which affects the $\mathfrak{R}$
operator. To show how an input operator influences the energy flows
we analyze a simplified version of magnetohydrodynamics, similar as
we did in \cite{SchoeberlICEAA2010}.

\subsection{Vibrating string}

Let us consider the system of a vibrating string. As independent coordinate
we choose the spatial coordinate $X$, the dependent coordinates are
the deflection $w$ and the temporal momentum $p$. This leads to
the following bundle $\mathcal{X}\rightarrow\mathcal{D},\,(X,w,p)\rightarrow X$,
The first jet manifold $\mathcal{J}^{1}(\mathcal{X})$ additionally
includes the derivative coordinates $w_{X}$ and $p_{X}$ and the
boundary $\partial\mathcal{D}$ consists of two points only, namely
$X=0$ and $X=L$ where $L$ is the length of the string. Approximately,
the pdes take the form 
\begin{equation}
\begin{array}{c}
\begin{array}{ccl}
\dot{w} & = & \frac{p}{\rho}\\
\dot{p} & = & d_{X}(P(X)w_{X})
\end{array}\end{array}\label{eq:PDEHeavyChain-1}
\end{equation}
with additional boundary conditions and the force in the string is
denoted by $P(X)$. To rewrite this system in a Hamiltonian fashion
we consider the Hamiltonian density $\mathcal{\mathfrak{H}=H}\mathrm{d}X$
with
\begin{equation}
\mathcal{H}=\frac{1}{2\rho}p^{2}+\frac{1}{2}P(X)w_{X}^{2}.\label{eq:HamChain-1}
\end{equation}
The total energy can be evaluated by $H=\intop_{0}^{L}\mathcal{H}\mathrm{d}X$.
To obtain partial differential equations in the form as in (\ref{eq:HamTimeInvStandardCF-1})
we can set $\mathcal{R}$ and $\mathbb{\mathcal{G}}$ to zero (no
damping and no inputs acting on the domain). Then we easily obtain
\[
\left[\begin{array}{c}
\dot{w}\\
\dot{p}
\end{array}\right]=\left[\begin{array}{cc}
0 & 1\\
-1 & 0
\end{array}\right]\left[\begin{array}{c}
\delta_{w}\mathcal{H}\\
\delta_{p}\mathcal{H}
\end{array}\right]=\mathcal{J}(\delta\mathfrak{H})
\]
with $x=(w,p)$ together with the expression 
\begin{equation}
\delta_{w}\mathcal{H}=\left(\partial_{w}-d_{X}\partial_{w}^{X}\right)\mathcal{H}=-d_{X}(P(X)w_{X})\label{eq:VarDerExamp}
\end{equation}
and 
\[
\delta_{p}\mathcal{H}=\left(\partial_{p}-d_{X}\partial_{p}^{X}\right)\mathcal{H}=\frac{p}{\rho}.
\]
Evaluating (\ref{eq:HdotPdeoO}) gives
\[
\dot{H}=\int_{\partial\mathcal{D}}\dot{x}\rfloor\delta^{\partial}\mathfrak{H}=\left.\dot{w}P(X)w_{X}\right|_{0}^{L}
\]
and again the power balance relation where the boundary conditions
are crucial in the determination of that expression.

Now we extend the equations by structural damping $r>0$ as
\begin{eqnarray*}
\dot{w} & = & \frac{p}{\rho}\\
\dot{p} & = & d_{X}(P(X)\, w_{X})+d_{X}(r\, d_{X}\frac{p}{\rho}).
\end{eqnarray*}
The port-Hamiltonian representation can be given as 
\[
\left[\begin{array}{c}
\dot{w}\\
\dot{p}
\end{array}\right]=\left(\left[\begin{array}{cc}
0 & 1\\
-1 & 0
\end{array}\right]-\left[\begin{array}{cc}
0 & 0\\
0 & -d_{X}\left(r\, d_{X}\left(\cdot\right)\right)
\end{array}\right]\right)\left[\begin{array}{c}
\delta_{w}\mathcal{H}\\
\delta_{p}\mathcal{H}
\end{array}\right]
\]
where now a differential operator is utilized in the $\mathfrak{R}$
mapping, as introduced in (\ref{eq:R_2ndorder_local}). The power
balance reads in this case as 
\begin{eqnarray*}
\dot{H} & = & -\int_{\mathcal{D}}\mathfrak{R}(\mathrm{\delta}\mathfrak{H)}\rfloor\delta\mathfrak{H}+\int_{\partial\mathcal{D}}\dot{x}\rfloor\delta^{\partial}\mathfrak{H}.
\end{eqnarray*}
Using (\ref{eq:Rintbyparts}) we obtain
\begin{eqnarray}
-\int_{\mathcal{D}}\mathfrak{R}(\mathrm{\delta}\mathfrak{H)}\rfloor\delta\mathfrak{H} & = & \int_{\mathcal{D}}\dot{w}d_{X}(r\, d_{X}\left(\dot{w}\right))\mathrm{dX}\label{eq:StringREq}
\end{eqnarray}
where the right hand side of (\ref{eq:StringREq}) decomposes into
\[
-\int_{\mathcal{D}}d_{X}\left(\dot{w}\right)r\, d_{X}\left(\dot{w}\right)\mathrm{dX}+\int_{\partial\mathcal{D}}\dot{w}\, r\, d_{X}\left(\dot{w}\right)\mathrm{d}X_{A}
\]
and thus 
\begin{eqnarray*}
\dot{H} & = & -\int_{\mathcal{D}}d_{X}\left(\dot{w}\right)\, r\, d_{X}\left(\dot{w}\right)\mathrm{dX}+\int_{\partial\mathcal{D}}\mathfrak{B}_{S}
\end{eqnarray*}
with the boundary expression 
\begin{eqnarray*}
\mathfrak{B}_{S} & = & \dot{w}\left(r\, d_{X}\left(\dot{w}\right)+P(X)\, w_{X}\right)\mathrm{d}X_{A}.
\end{eqnarray*}
It is obvious that the structural damping extracts power as can be
seen from the term on the domain
\[
-\int_{\mathcal{D}}d_{X}\left(\dot{w}\right)\, r\, d_{X}\left(\dot{w}\right)\mathrm{dX}=-r\int_{\mathcal{D}}\left(d_{X}\left(\dot{w}\right)\right)^{2}\mathrm{dX}
\]
in $\dot{H}$ since $r>0$.

\subsection{Magnetohydrodynamics (MHD)}

In this section we consider an example that enables us to introduce
a differential operator in the in/output mapping, namely so-called
inductionless MHD. Roughly speaking, we analyze the macroscopic behavior
of an electrically conducting fluid in the presence of external electromagnetic
fields. We assume that the dynamic of the additionally induced electromagnetic
parts can be neglected (low magnetic Reynold\textquoteright{}s number).
Furthermore, we confine ourselves to the case of negligible viscosity
and electrical conductivity. 

To obtain the mathematical model we consider again a bundle structure
of the form $\mathcal{X}\rightarrow\mathcal{D}\,,\,(X^{A},q^{\alpha})\rightarrow(X^{A})$
and we introduce the so-called material mass and charge density, $\rho\left(X\right)$
and $\mu\left(X\right)$. Furthermore we assume the existence of a
stored energy function $E_{st}\left(q_{A}^{\alpha}\right)$ which
meets 
\begin{eqnarray*}
\mathcal{P} & = & -\rho\frac{\partial E_{st}}{\partial J}\,,\,\,\,\,\,\, J=\det\left(F\right)
\end{eqnarray*}
with the deformation gradient $F_{A}^{\alpha}=q_{A}^{\alpha}$, and
$F_{A}^{\alpha}\hat{F}_{\alpha}^{B}=\delta_{A}^{B}$, where we have
introduced the material pressure $\mathcal{P}\left(X\right)$. The
electrostatic potential is denoted by $A_{0}\left(q\right)$ and the
vector potential by $A_{\alpha}\left(q\right)$, where the components
of the electric field strength are given by $E_{0\alpha}=\partial_{\alpha}A_{0}$
and of the magnetic flux density by 
\begin{eqnarray*}
B_{\alpha\beta} & = & \partial_{\alpha}A_{\beta}-\partial_{\beta}A_{\alpha}.
\end{eqnarray*}
The Hamiltonian density $\mathfrak{H}=\mathcal{H}\mathrm{d}X$ follows
as 
\begin{eqnarray*}
\mathcal{H} & = & \frac{1}{2\rho}\delta^{\alpha\beta}(p_{\alpha}-\mu A_{\alpha})(p_{\beta}-\mu A_{\beta})+\rho E_{st}
\end{eqnarray*}
with the momenta 
\begin{eqnarray*}
p_{\alpha} & = & \rho\delta_{\alpha\beta}\dot{q}^{\beta}+\mu A_{\alpha}.
\end{eqnarray*}

\begin{rem}
The Hamiltonian can be derived using the Lagrangian 
\[
\mathcal{L}=\frac{1}{2}\rho\dot{q}^{\alpha}\delta_{\alpha\beta}\dot{q}^{\beta}-\rho E_{st}+\mu\left(A_{0}+\dot{q}^{\alpha}A_{\alpha}\right)
\]
together with $\mathcal{H}=\left(p_{\alpha}\dot{q}^{\alpha}-\mathcal{L}\right)$
where the potential $A_{0}$ is neglected in a first step, since it
serves as the input, which is considered separately. 
\end{rem}
The partial differential equations read as
\[
\dot{p}_{\alpha}=-J\hat{F}_{\alpha}^{A}\partial_{A}\mathcal{P}+\mu\frac{1}{\rho}\delta^{\beta\rho}(p_{\rho}-\mu A_{\rho})\partial_{\alpha}A_{\beta}
\]
and its Hamiltonian representation follows as
\[
\left[\begin{array}{c}
\dot{q}^{\alpha}\\
\dot{p}_{\alpha}
\end{array}\right]=\left[\begin{array}{cc}
0 & \delta_{\beta}^{\alpha}\\
-\delta_{\alpha}^{\beta} & 0
\end{array}\right]\left[\begin{array}{c}
\delta_{\beta}\mathcal{H}\\
\delta^{\beta}\mathcal{H}
\end{array}\right].
\]
with $\delta^{\beta}\mathcal{H}=\partial^{\beta}\mathcal{H}=\frac{\partial}{\partial p_{\beta}}\mathcal{H}$. 

Taking the system input $A_{0}$ into account we extend the system
and obtain the following representation 
\begin{eqnarray*}
\left[\begin{array}{c}
\dot{q}^{\alpha}\\
\dot{p}_{\alpha}
\end{array}\right] & = & \mathcal{J}\left[\begin{array}{c}
\delta_{\beta}\mathcal{H}\\
\delta^{\beta}\mathcal{H}
\end{array}\right]+\mathfrak{G}\left(A_{0}\right)\,,
\end{eqnarray*}
and 
\[
\mathfrak{G}\left(A_{0}\right)=\left[\begin{array}{c}
0\\
\mu\hat{F}_{\alpha}^{B}d_{B}A_{0}
\end{array}\right],
\]
where the input map corresponds to the definition of equation (\ref{eq:G_1storder_local})
with $u=A_{0}$. 

The power balance then again takes the form
\[
\dot{H}=\int_{\mathcal{D}}\mathfrak{G}(A_{0})\rfloor\delta\mathfrak{H}+\int_{\partial\mathcal{D}}\dot{x}\rfloor\delta^{\partial}\mathfrak{H}
\]
where the interesting term corresponds to the pairing on the domain
including the input and we obtain 
\begin{eqnarray*}
\int_{\mathcal{D}}\mathfrak{G}(A_{0})\rfloor\delta\mathfrak{H} & = & \int_{\mathcal{D}}\mu\hat{F}_{\alpha}^{B}(d_{B}A_{0})\dot{q}^{\alpha}\mathrm{d}X\\
 & = & -\int_{\mathcal{D}}d_{B}(\hat{F}_{\alpha}^{B}S^{\alpha})A_{0}\mathrm{d}X+\int_{\partial\mathcal{D}}\mathfrak{B}_{M}
\end{eqnarray*}
with 
\begin{eqnarray*}
\mathfrak{B}_{M} & = & S^{\alpha}\hat{F}_{\alpha}^{B}A_{0}\mathrm{d}X_{B},
\end{eqnarray*}
see also (\ref{eq:Gintbyparts}). Here the adjoint operator on the
domain takes the form $\mathfrak{G}^{*}\left(\delta\mathfrak{H}\right)=-d_{B}(\hat{F}_{\alpha}^{B}S^{\alpha})$.

The expression 
\begin{eqnarray*}
S^{\alpha} & = & \frac{\mu}{\rho}\delta^{\alpha\beta}\left(p_{\beta}-\mu A_{\beta}\right)
\end{eqnarray*}
 is equivalent to the components of the convective current density
which is caused by the convective transport of charge.

\section{Conclusion}

Based on the relation (\ref{eq:LieDiffPde}) that describes the power
balance relation, we introduced port-Hamiltonian systems modeled by
partial differential equations, in such a way that energy conservation,
dissipation, collocation and energy ports come apparent. This key
idea is very similar of that used in the definition of the Stokes-Dirac
structure, however due to a different choice of coordinates (not necessarily
energy coordinates) the port-Hamiltonian representations of the two
approaches differ significantly. This can be observed by a close inspection
of the differential operators that are needed to describe the physical
systems properly. Therefore it is of great interest for future investigations
two analyze and compare these two promising directions more closely.

\section*{Appendix}

Les us consider a bundle $\pi:\mathcal{X}\rightarrow\mathcal{D},\,(X^{A},x^{\alpha})\rightarrow(X^{A})$
and its n-th order jet manifold $\mathcal{J}^{n}(\mathcal{X})$ as
well as the bundle $\pi_{0}^{n}:\mathcal{J}^{n}(\mathcal{X})\rightarrow\mathcal{X}$
together with a section $\phi:\mathcal{D}\rightarrow\mathcal{X}$,
i.e. $x=\phi(X)$. 

The exterior derivative $\mathrm{d}$ is connected to the horizontal
derivative $\mathrm{d}_{h}$ through
\begin{equation}
\left(j^{n+1}\phi\right)^{\ast}\left(\mathrm{d}_{h}\left(\omega\right)\right)=\mathrm{d}\left(\left(j^{n}\phi\right)^{\ast}\left(\omega\right)\right)\label{eqHorizontDeriv}
\end{equation}
for a form $\omega$ living on $\mathcal{J}^{n}(\mathcal{X})$. Roughly
speeking the pull back $\left(j^{n}\phi\right)^{\ast}\left(\omega\right)$
denotes $\omega\circ(j^{n}\phi)$ where $j^{n}\phi$ corresponds to
the n-th order jet-lift of the section $\phi$. 

Furthermore we have 
\[
\int_{\mathcal{D}}j^{n+1}\left(\phi\right)^{\ast}\left(\mathrm{d}_{h}\omega\right)=\int_{\mathcal{D}}\mathrm{d}\left(j^{n}\left(\phi\right)^{\ast}\left(\omega\right)\right)=\int_{\partial\mathcal{D}}j^{n}\left(\phi\right)^{\ast}\left(\omega\right)
\]
for $\omega\in(\pi_{0}^{n})^{*}(\Lambda^{d-1}(\mathcal{X}))$, which
is nothing else than Stokes' theorem adapted to bundles.

\begin{ack}
Markus Sch\"oberl is an APART fellowship holder of the Austrian Academy of Sciences.
\end{ack}

\bibliography{maxibib}             
                                                
\end{document}